\documentclass[10pt]{amsart}

\usepackage{amssymb,eucal}
\usepackage{amscd}
\usepackage[all,cmtip]{xy}
\usepackage{rotating}
\usepackage{hyperref,mathrsfs,soul}
\usepackage{tikz}

\usepackage{amsmath,amssymb}

\setul{}{1.3pt}

\newtheorem{theorem}{Theorem }[section]

\newtheorem{lemma}[theorem]{Lemma}

\newtheorem{remark}[theorem]{Remark}
\newtheorem{corollary}[theorem]{Corollary}

\newcommand{\proj}{\mathrm{proj}}
\def\I{\mathbf{I}}

\def\eop{\hspace*{\fill}{\footnotesize$\blacksquare$}}
\newcommand{\Aut}{\mathrm{Aut}}
\newcommand{\id}{\mathrm{id}}

\newcommand{\mI}{\mathcal{I}}

\newcommand{\mP}{\mathcal{P}}

\newcommand{\mM}{\mathcal{M}}
\newcommand{\mS}{\mathcal{S}}

\newcommand{\mF}{\mathcal{F}}

\newcommand{\mB}{\mathcal{B}}
\newcommand{\bP}{\mathbb{P}}
\newcommand{\mL}{\mathcal{L}}
\newcommand{\mW}{\mathcal{W}}
\newcommand{\F}{\mathbb{F}}
\newcommand{\mH}{\mathcal{H}}

\title[A question of Frohardt]{A question of Frohardt on $2$-groups, skew translation quadrangles of even order and cyclic STGQs}


\author{Koen Thas}

\address{{Ghent University},
{Department of Mathematics},
{Krijgslaan 281, S25, B-9000 Ghent, Belgium}}
\email{koen.thas@gmail.com}
\thanks{}
\date{}

\begin{document}
\maketitle

\begin{abstract}
 We solve a fundamental question posed in Frohardt's 1988 paper \cite{Fro} on finite $2$-groups with Kantor familes, by showing that finite groups $K$ with 
 a Kantor family $(\mF,\mF^*)$ having distinct members
$A, B \in \mF$ such that $A^* \cap B^*$ is a central subgroup of $K$ and the quotient $K/(A^* \cap B^*)$ is abelian cannot exist if the center of $K$ has exponent $4$ and the members of $\mF$ are elementary abelian. Then we give a short geometrical proof of a recent result of Ott  
which says that finite skew translation quadrangles of even order $(t,t)$ (where $t$ is not a square) are always translation generalized quadrangles. This is a consequence of a complete classification of finite cyclic skew translation quadrangles of order $(t,t)$ that we carry out in the present paper. 
\end{abstract}

\setcounter{tocdepth}{1}
\tableofcontents

\bigskip
\section{Introduction}

In a seminal paper of 1988 \cite{Fro}, Frohardt solved a big part of a famous question of Kantor which asks whether a finite group admitting 
a Kantor family is a $p$-group. Such groups are very important since they produce generalized quadrangles, the central Tits buildings 
of rank $2$. In fact, groups with Kantor families yield the most important tool to construct generalized quadrangles, and a large 
literature exists on these objects | see for instance \cite{PT,TGQBook,LEGQ,HB,TiWe}.  Generalized quadrangles constructed from Kantor families | by a procedure which is 
explained in section \ref{KF} | carry an interesting automorphism group which locally fixes some point linewise, and are called ``elation generalized quadrangles'' (EGQs). 
The theory of EGQs was initially devised through the 1970s and 1980s (see \cite[chapters 8--10]{PT} and \cite{LEGQ}) in order to set up a framework 
to study automorphisms, and especially Moufang conditions, of generalized quadrangles from a local point of view, much like translation planes were studied 
in the theory of projective planes. Together with a deep combinatorial theory, this enabled an almost entirely synthetic-geometric proof of the celebrated results 
of Fong and Seitz \cite{BN1,BN2} on split BN-pairs of rank $2$, in the case of type $\texttt{B}_2$ BN-pairs | see \cite[chapters 8 and 9]{PT}, \cite[chapters 10 and 11]{TGQBook}, and also the related recent note \cite{localhalf}. \\

In \cite{Fro}, Frohardt also studied groups $K$ with Kantor families $(\mF,\mF^*)$ which have distinct members
$A, B \in \mF$ such that $A^* \cap B^* \leq Z(K)$ (where $Z(K)$ denotes the center of $K$) and $K/(A^* \cap B^*)$ is abelian. Such Kantor families are fundamental, since almost 
all examples of Kantor families have this natural property | see, e.g., the monograph \cite{LEGQ} or Payne's famous paper \cite{9}. We also refer to 
section \ref{STGQover} for more background information on this matter. Only a very few classes of finite groups are known which have Kantor families; the most used 
ones are elementary abelian groups, and Heisenberg groups defined over finite fields in odd characteristic, of dimension $3$ or $5$ (see e.g. the monograph \cite{LEGQ}). They all satisfy Frohardt's condition.

\begin{theorem}[D. Frohardt \cite{Fro}, 1988]
\label{Fro88}
Let $K$ be a finite group with a Kantor family $(\mF,\mF^*)$ having distinct members
$A, B \in \mF$ such that $A^* \cap B^* \leq Z(K)$ and $K/(A^* \cap B^*)$ is abelian. Then one of the folowing three cases occurs:
\begin{itemize}
\item[(1)]
$Z(K)$ and the elements of $\mF$ are elementary abelian.
\item[(2)]
$Z(K)$ is an elementary abelian $2$-group and the elements of $\mF$ have exponent $4$.
\item[(3)]
$Z(K)$ has exponent $4$ and the elements of $\mF$ are elementary abelian $2$-groups.
\end{itemize}
\end{theorem}

In {\em loc. cit.}, Frohardt asked whether cases (2) and (3) actually occur. \\

In 2006, Rostermundt \cite{Roster} and, independently Thas \cite{Basic}, constructed an infinite class of examples  (related to Hermitian varieties in
projective dimension $3$), which fitted in class (2), so class (3) became the final challenge.
Frohardt's problem resisted many years of tries (see for instance the work of Hachenberger \cite{H}), and the problem also occurs in the large literature of translation nets. Besides, the existence of elements in class (3) presented one of the main obstacles in classifying so-called ``skew translation generalized quadrangles'' (STGQs) | see \cite{noteann,STGQ}, and also \cite{STGQ2,Leug}. We also refer to section \ref{STGQover} for some milestones in that classification theory. \\

  


In this paper, we resolve the question by showing that class (3) is empty. 
The proof consists of a mixture of synthetic and geometric reasoning. 
Payne's first paper on skew translation generalized quadrangles \cite{STGQbirth} from 1975, marking the birth of the theory.  In its long history, no examples 
were found of STGQs of even order $(t,t)$ for which the elation group was {\em not} abelian | when $t$ is odd, the theory is entirely different, since then it is easy to prove that 
such a group cannot be abelian at all. Since one expects that such STGQs would have to meet the requirements of Frohardt's theorem, and since for STGQs of order $(t,t)$
the elements of $\mF$ are always elementary abelian (see \cite{STGQ} for those details),  such examples would conjecturally live in class (3) of Theorem \ref{Fro88}.
And indeed they do, if they would exist.
The best result up till 2020 was the same result when $(t,t) = (8,8)$ \cite{Leug}, with a very long and computer-aided technical proof.  In the meanwhile, Ott has published a long paper \cite{Ott} in which he showed that finite skew translation quadrangles of order $(t,t)$ with $t$ not a square, are translation generalized quadrangles, that is, the associated elation group is an elementary 
abelian $2$-group, and hence the quadrangle has a projective representation in some projective space. This essentially means that there is no ``proper'' theory of 
STGQs in this case. The present paper contains an alternative approach to the one of Ott, although we use one lemma of his paper. 
Our proof, which is geometric in nature, is very short in comparison with Ott's proof, and we first completely classify skew translation generalized quadrangles of order $(t,t)$ which 
have a cyclic group of symmetries with center the elation point. We then use the lemma of Ott to reduce the general problem to the case 
of ``cyclic STGQs.''

\begin{remark}{\rm 
On the combinatorial level, there is also an other interesting aspect to this result. If every line incident with some point in a finite generalized quadrangle of even order $(t,t)$ 
is {\em regular} (see section \ref{reg} for a formal definition of this very important notion), then it can be shown that the point itself is also regular (\cite[section 1.5.2]{PT}). Ott's result provides a group-theoretical converse if $t$ is not a square: if a finite EGQ
of even order $(t,t)$ with $t$ not a square, has a regular elation point (this is another way of defining STGQs in this specific case), then it is a translation quadrangle, so all the lines incident with $x$ are regular.  }
\end{remark}

\subsection*{Organization}

In sections \ref{intro} and \ref{sett}, we introduce the basic notions that we will need. In  section \ref{prep}, we make a number of synthetic observations; 
we essentially show that in case (3) of Frohardt's theorem, it is sufficient to work with Kantor families of type $(t,t)$. This is done in the geometrical language of 
generalized quadrangles. Then, in section \ref{Frohthm}, we show that Frohardt's class (3) is empty. 
And finally, in the last section, we classify cyclic STGQs (of finite order $(t,t)$). 
As a consequence, STGQs of even order $(t,t)$ with $t$ not a square are always translation quadrangles.

\bigskip
\section{Synopsis of definitions}
\label{intro}

Let $\Gamma$ be a thick generalized quadrangle (GQ). It is a rank $2$ geometry $\Gamma = (\mP,\mL,\I)$ (where we call the elements of $\mP$ ``points'' and those of $\mL$ ``lines'')
such that the following axioms are satisfied:
\begin{itemize}
\item[(a)]
there are no ordinary digons and  triangles contained in $\Gamma$;
\item[(b)]
each two elements of $\mP \cup \mL$ are contained in an ordinary quadrangle;
\item[(c)]
there exists an ordinary pentagon.
\end{itemize}
It can be shown that there are exist constants $s$ and $t$ such that each point is incident with $t + 1$ lines and each line is incident with $s + 1$ points. We say that
$(s,t)$ is the {\em order} of $\Gamma$. 
Note that an ordinary quadrangle is just a (necessarily thin) GQ of order $(1,1)$ | we call such a subgeometry also ``apartment'' (of $\Gamma$).

\subsection{Subquadrangles}

A {\em subquadrangle} (subGQ) $\Gamma' = (\mP',\mL',\I')$ of a generalized quadrangle $\Gamma = (\mP,\mL,\I)$, is a GQ for which 
$\mP' \subseteq \mP$, $\mL' \subseteq \mL$, and $\I'$ is the incidence relation which is induced by $\I$ on 
$(\mP' \times \mL') \cup (\mL' \times \mP')$.

\subsection{Regularity}
\label{reg}

Let $\Gamma$ be a thick GQ of order $(s,t)$. If $X$ and $Y$ are (not necessarily different) points, or lines, $X \sim Y$ denotes 
the fact that there is a line incident with both, or a point incident with both. Then $X^{\perp} := \{ Y \ \vert\ Y \sim X \}$, and 
if $S$ is a point set, or a line set, $S^{\perp} := \cap_{s \in S}s^{\perp}$ and $S^{\perp\perp} := {(S^{\perp})}^{\perp}$. 

A particularly important example is the case where $S = \{X,Y\}$ is a set of distinct noncollinear points (or nonconcurrent lines, but this is 
merely the dual situation, which we leave to the reader); then each line incident with $X$ is incident with precisely one point of $\{ X,Y\}^{\perp}$ (so if $\Gamma$ is finite, 
$\vert \{ X,Y\}^{\perp} \vert = t + 1$). The set $\{ X,Y\}^{\perp\perp}$ consists of all points which are collinear with every point of $\{X,Y\}^{\perp}$, 
so 
\begin{equation}
\{ X, Y\} \subseteq \{X,Y\}^{\perp\perp}. 
\end{equation}

Let $Z$ be any point of $\{X,Y\}^{\perp}$; if each line incident with $Z$ is incident with exactly one point of $\{X,Y\}^{\perp\perp}$, then 
this property is independent of the choice of $Z$, and we say that $\{ X,Y\}$ is {\em regular}. In the finite case, we could equivalently have asked 
that $\vert \{ X,Y\}^{\perp\perp} \vert = t + 1$. 

We call a point/line $X$ {\em regular} if $\{ X,Y \}$ is regular for all $Y \not\sim X$.

\subsection{Symmetry}

Isomorphisms and automorphisms of generalized quadrangles are defined in the usual manner. See chapter 1 of \cite{PT}. The automorphism group of a 
GQ $\Gamma$ will be denoted by $\Aut(\Gamma)$. 

Let $X$ be a point or a line in a thick GQ $\Gamma$. A {\em symmetry} with {\em center} $X$ (in the case of a point) or {\em axis} $X$ (in the case of a line)
is an element of $\Aut(\Gamma)$ that fixes each element of $X^{\perp}$. We say that $X$ is a {\em center of symmetry} (point case) or an {\em axis of symmetry} (line 
case) if there exist  $Y$ and $Z$  in $X^{\perp}$ such that $Y \not\sim Z$, for which the group of all symmetries $\mS(X)$ with center/axis $X$ acts 
transitively on $\{ Y,Z\}^{\perp} \setminus \{X\}$. This definition does not depend on the choice of $(Y,Z)$, and one easily shows that ``transitive'' implies
``sharply transitive.'' In the finite, we could also have asked that 
\begin{equation}
\vert \mS(X) \vert = t
\end{equation}
if the order of $\Gamma$ be $(s,t)$.

Note that a center/axis of symmetry is necessarily regular.\\

We end this section with a very useful result of Benson. 

\begin{theorem}[C. T. Benson \cite{benson}]
\label{benson}
Let $\Gamma$ be a thick finite GQ of order $(s,t)$, and let $\gamma$ be an automorphism of $\Gamma$. If $N$ is the number of fixed points of $\gamma$, and $M$ is the number of points $x$ for which $x^\gamma \ne x \sim x^\gamma$, then
\begin{equation}
(t + 1)N + M \equiv st + 1 \mod{(s + t)}.
\end{equation}
\end{theorem}

\subsection{Kantor families}
\label{KF}

In this section, we will recall the notion of Kantor family construction. Groups with Kantor families produce generalized quadrangles, and 
vice versa, a certain type of generalized quadrangle gives rise to Kantor families. We will only define {\em finite} Kantor families, as that is the case 
we will need in due course. 

So suppose $K$ is a finite group of order $u^2v$ for positive integers $u$ and $v$, both at least $2$. A {\em Kantor family} (of type $(u,v)$) 
in $K$ is a pair $(\mF,\mF^*)$ of sets of subgroups of $K$ for which the properties below are satisfied:
\begin{itemize}
\item[(a)]
$\vert \mF \vert = \vert \mF^* \vert = v + 1$, and there is a bijection $*: \mF \mapsto \mF^*$ which maps each $A \in \mF$ to 
$A^* \in \mF^*$ such that $A \leq A^*$;
\item[(b)]
for each $A \in \mF$, $\vert A \vert = u$ and $\vert A^* \vert = uv$;
\item[(c)]
if $A, B, C$ are different elements in $\mF$, then $AB \cap C = \{ \id \}$;
\item[(d)]
if $A$ and $B$ are different elements in $\mF$, then $A \cap B^* = \{ \id \}$. 
\end{itemize}

From the data $\Big(K,(\mF,\mF^*)\Big)$ one constructs a thick GQ $\Gamma = \Gamma\Big(K,(\mF,\mF^*)\Big)$ of order $(u,v)$ as follows. Its points are a symbol $(\infty)$, left cosets of type $gA^*$ ($A \in \mF$ and $g \in K$), 
and the elements of $K$. Its lines are symbols $[ A ]$ ($A \in \mF$) and left cosets of type $gA$ ($A \in \mF$ and $g \in K$). The point $(\infty)$ is incident 
with all lines of the first type, and no other lines. A line $[A]$ is also incident with all points $gA^*$. 
All other incidences are (reversed) containment. Elements of the group $K$ act naturally as automorphisms  
of $\Gamma$, by left multiplication on the cosets, while fixing the symbolic point and symbolic lines. Then $K$ fixes all the lines incident with $(\infty)$, 
and acts sharply transitively on the points of $\Gamma$ which are not collinear with $(\infty)$. 

Conversely, let $\Omega$ be a thick GQ of order $(u,v)$, with an automorphism group $L$ which fixes some point $x$ linewise while acting sharply transitively 
on the points not collinear with $x$. We call such GQs  {\em elation generalized quadrangles} (EGQs) with {\em elation point} $x$ and {\em elation group} $L$. 
Then $\vert L \vert = u^2v$. Now take one arbitrary point $y$ not collinear with $x$ (due to the transitivity of $L$ on these 
points, this choice is not important). For each line $U \I y$, let $\widetilde{u}$ be the point which is incident with $U$ and collinear with $x$; let $L_U^* := L_{\widetilde{u}}$.
then with 
$\mF := \{ L_U\ \vert\ U \I y \}$ and $\mF^* := \{ L_U^*\ \vert\ U \I y \}$, $(\mF,\mF^*)$ is a Kantor family of type $(u,v)$ in $L$. Also, we have a natural isomorphism
\begin{equation}
\Gamma\Big(L,(\mF,\mF^*)\Big)\ \mapsto\ \Omega
\end{equation}
which maps $L$ to itself and $(\infty)$ to $x$.

\subsection{Skew translation quadrangles}
\label{STGQover}

If $(\Gamma,K)$ is an EGQ with elation point $x$, and $x$ is a center of symmetry such that the corresponding symmetry group $\mS$ is a subgroup of $K$, then 
we call $(\Gamma,K)$ a {\em skew translation generalized quadrangle} (STGQ). In terms of Kantor families (using the notation of the previous subsection), a Kantor family 
$(\mF,\mF^*)$ gives rise to an STGQ if and only if there is a normal subgroup $C$ of $K$ such that for each $A \in \mF$, we have $A^* = AC$ (see \cite{STGQbirth} or \cite[chapter 10]{PT});  $C$ then corresponds to the group of symmetries with center the elation point. Note that for different $A, B$ in $\mF$, we have $A^* \cap B^* = C$. This type of GQ is very general: each known finite generalized quadrangle which is not isomorphic to the Hermitian quadrangle $\mH(4,q^2)$ for some finite field $\F_q$, 
is | modulo point-line duality | either an STGQ, or the Payne-derived quadrangle of an STGQ. More details on this statement can be found in \cite{STGQ}; for a more or less recent but rather detailed census 
on the known GQs, we also refer to \cite{Payneproc}. 

We recall some basic (``decisive'') classification results in the large theory of STGQs. Much more can be found in \cite{STGQ}. The first result was first obtained by the author in \cite{STGQ} (see \cite{noteann} for chronological details in its conception). It was also obtained in \cite{Leug}
with virtually the same proof. 

\begin{theorem}
\label{oddsq}
An STGQ of odd order $(t,t)$ is isomorphic to the symplectic quadrangle $\mW(t)$. 
\end{theorem}

The symplectic quadrangle $\mW(t)$ arises as the $\F_t$-rational points in a projective space $\bP^3(\F_t)$, together with the absolute lines with respect to 
a symplectic polarity of $\bP^3(\F_t)$ (see \cite[chapter 3]{PT}); it is one of the natural geometric modules of the group $\mathbf{PSp}_4(t)$.\\

The most investigated type of STGQ is arguably the class of ``flock quadrangles.'' All the known flock quadrangles arise through Kantor families in one and the same 
type of elation group: finite Heisenberg groups of dimension $5$ (over a finite field). The next result proves the converse, and is taken from \cite{Isomflock}.

\begin{theorem}
If $(\Gamma,K)$ is an EGQ of order $(s,t)$ and $K$ is isomorphic to a Heisenberg group of dimension $5$ over $\F_t$, then $\Gamma$ is a flock quadrangle. 
\end{theorem}

The following result is one of the only known results on STGQs of even order $(t,t)$; its proof takes up more than half of the paper \cite{Leug}, and 
is computer-aided. 

\begin{theorem}
An STGQ of order $(8,8)$ is a translation generalized quadrangle.
\end{theorem}

Finally, we have Ott's beautiful result:

\begin{theorem}[Ott \cite{Ott}]
An STGQ $(\Gamma,K)$ of even order $(t,t)$ is a TGQ if $t$ is not a square. In that case, $K$ is an elementary abelian $2$-group.
\end{theorem}


\medskip
\section{Notation}
\label{sett}

Let $K$ a finite group with a Kantor family $(\mF,\mF^*)$ of type $(s,t)$. The associated generalized quadrangle with parameters $(s,t)$ is denoted $\Gamma^x$, or $\Gamma$, where the elation point is $x$. The Kantor family is defined relative to 
the point $z \not\sim x$ (unless we mention it otherwise).

We say that a group $G \leq K$ satisfies (*) if for any element $\alpha$ of $G$, we have that if $\alpha$ fixes some 
point $y \sim x \ne y$, $\alpha$ also fixes $xy$ pointwise. We say that $G \leq K$ satisfies (*) {\em at} $U \I x$
if this property is locally fulfilled at $U$.

In section \ref{prep} and section \ref{Frohthm}, we assume that for certain distinct elements $A, B \in \mF$, we suppose that $K/(A^* \cap B^*)$ is abelian and that $A^* \cap B^* \leq Z(K)$, the center of $K$. Put $A^* \cap B^* =: \mS$. 
We want to work solely in case (3) of Frohardt's theorem in sections \ref{prep} and \ref{Frohthm}, so that the elements of $\mF$ are elementary abelian, and  
$Z(K) =: Z$ has exponent $4$. This means that $\Gamma$ is not a TGQ (and hence $K$ is not abelian).

\medskip
\section{Preliminary properties, and reduction to central STGQs of even order $t$} 
\label{prep}

In this section we observe some initial properties which narrow down the possibilities for elements in Theorem \ref{Fro88}(3). In the final part of the proof of the main 
theorem, we will not need all these properties, but it is interesting to see how far one can go ``synthetically.'' \\

First note that for any line $W \I x$, $K$ acts transitively on $W \setminus \{x\}$. As $\mS \leq Z(K)$, it hence follows that 
$\mS$ fixes $[A]$ and $[B]$ pointwise (so $\mS$ is in the kernel of the action of $K$ on the point set of $[A]$, and $[B]$). As $K/\mS$ is abelian
and transitive on both $[A] \setminus \{x\}$ and $[B] \setminus \{x\}$, we have that $K$ has (*) at $[A]$ and $[B]$.

Let $z \in Z^\times$; if $z$ would fix an affine line (which is by definition a line not incident with $x$), $z$ is in some conjugate of an element of $\mF$, implying that $z^2 = 1$. 
Now we assume that $z$ does not have this property. Now suppose there is an affine line $U$ such that $U^z$ does not intersect $U$.
Take any element $V \in \{U,U^z\}^{\perp}$ which is not incident with $x$; choose the element $\beta$ in $K_V$ which sends $U^z$ to $U$; then 
$z\beta$ fixes $U$, and $[z,\beta] = \id$, so $\id = (z\beta)^2 = z^2\beta^2 = z^2$. So we may suppose that $z$ fixes each point of 
$x^{\perp}$, i.e., $z$ is a symmetry with center $x$.  In its turn, this implies that $z \in \mS$. \\

\subsection{Reduction of parameters}

Now suppose that $z$ is a (nontrivial) symmetry with center $x$, and suppose that $U \sim [A]$ is an affine line. 
Let $\alpha \in A^*$ be an element which sends $U^z$ to $U$, such that $z\alpha \ne \id$ (and note that such $\alpha$ exist, of course).
Then $z\alpha$ fixes $U$. As $[z,\alpha] = \id$, $z$ is an involution if and only if $\alpha$ is an involution. As $(\mF,\mF^*)$ is a Kantor family,
we can write $A^* = A\mS$ (with $A \cap \mS = \{\id\}$). We can write $\alpha = z^{-1}a$, with $a \in K_U$, so that the fixed points structure
of $\alpha$ is that of $a$ | precisely the points incident with $[A]$. Suppose $\alpha$ does not fix affine lines.
Applying Theorem \ref{benson}, we conclude
\begin{equation}
(t + 1)(s + 1) + st\ \equiv\ st + 1 \mod{(s + t)},
\end{equation}
which gives a contradiction unless $s = t$ is even.
In the other cases, $\alpha$ must fix affine lines, so that it is an involution, and hence $z$ also is. So if $s \ne t$, 
each element of $Z$ is an involution, showing that case (3) of Frohardt's theorem cannot occur under these assumptions. {\bf it follows that 
$s = t$ is even.}

\subsection{The centralizers $C_B(\alpha)$}
\label{cent}

Now consider the case where $s = t$, and $t$ is even. Note that $\mS$ only consists of symmetries and involutions. 
Suppose $\beta \in \mS^\times$ is not an involution (it is then a symmetry); compose $\beta$ with a hypothetical involution $\iota$ in $Z$ which is not a symmetry.
Then $\beta\iota \in Z$ is not a symmetry, so it is an involution, and $\id = (\beta\iota)^2 = \beta^2\iota^2 = \beta^2$, contradiction. So all 
involutions in $\mS$ are symmetries, and whence \ul{all elements of $\mS$ are symmetries.} It follows easily that $x$ is a regular point.  

If $K$ is abelian, the theory of translation generalized quadrangles gives us that $K$ is elementary abelian (see \cite[chapter 3]{TGQBook}), so we may assume that $K$ is not 
abelian. {\bf By Proposition 3.1 of Hachenberger \cite{H}, it follows that $Z(K) = \mS$.}  \\

Let $C \ne D$ be in $\mF$; then $[K,K] = [C,D]$ (as $K = CD\mS$) $= \langle (cd)^2 \vert c \in C, d \in D \rangle$. As this group is 
a subgroup of $\mS$, it follows that $[K,K]$ is elementary abelian. As $\Phi(K) = K^2[K,K]$ (where $\Phi(K)$ denotes the Frattini subgroup of $K$), it easily follows that $\Phi(K)$ is the elementary abelian subgroup of $\mS$
of all squares. 
Then on the other hand, if $s = \gamma^2 \in \mS$ is a square, there is an $e \in E \in \mF$ and $c \in \mS$ such that $\gamma = ec$ (as 
$K = \cup_{A \in \mF}A^*$), and $s = ecec = c^2$. So $\Phi(K) = \mS^2$. Also, for $A \in \mF \setminus \{ E\}$, there is a unique $B \in \mF$ such that 
$\gamma \in AE$ | write $\gamma   ab$ with $a \in A$ and $b \in B$. Then $(ab)^2 = \gamma^2 = s$, so that $s \in [K,K]$. {\bf We have obtained 
that $[K,K] = \mS^2 = \Phi(K)$.}\\

Let $A \ne B$ be arbitrary in $\mF$.
Now let $\alpha \in A^\times$ be arbitrary, and let $\beta \in B^\times$ be such that $[\alpha,\beta] = \id$. Put $\mF \setminus \{A,B\} = \{C_1,\ldots,C_{t - 1}\}$.
Then for each $i \in \{1,\ldots,t - 1\}$ there is precisely one triple $(c_i,\beta_i,s_i) \in C_i \times B \times \mS$ such that $\alpha\beta = 
c_i\beta_is_i$. Note that the maps $\mu: \{ 1,\ldots,t - 1\} \longrightarrow \mS^\times: j \longrightarrow s_j$ and $\mu': \{ 1,\ldots,t - 1\} \longrightarrow B \setminus \{\beta\}: j \longrightarrow \beta_j$
are surjective. 
Suppose that $s_e^2 = \id$; then as $(\alpha\beta)^2 = \id$, $[c_e,\beta_e] = \id$, and so $[\beta_e,c_e\beta_e s_e] = \id = [\beta_e,\alpha\beta]$, 
implying that $[\alpha,\beta_e] = \id$. The converse is also true. So 
\begin{equation}
\vert \mbox{involutions}\ \mbox{in}\ \mS\vert (\mbox{including} \ \id) = \vert C_B(\alpha)\vert =: \ell.
\end{equation}

On the other hand we have $t = \vert \alpha^B \vert \times \vert C_B(\alpha) \vert$, so for any such $\alpha$, the size 
of $\alpha^K$ (which is the same as $\alpha^B$) is a constant ($\frac{t}{\ell}$) independent of the choice of $\alpha$. As $\alpha^K = \alpha^B$, 
it also follows that $\alpha$ fixes precisely $\vert C_B(\alpha) \vert = \ell$ points of $[A] \setminus \{x\}$ linewise.

Note at this point that if we prove that $\mS$ is an elementary abelian $2$-group, then each $A^*$ is also elementary abelian, and as $K$ is covered by the $A^*$s, it follows that $K$ itself is elementary abelian, contradiction with the assumption that $K$ has exponent $4$.

In the rest of this paper, we will denote the subgroup of $\mS$ that consists of all its involutions, by $\mI$.

\subsection{Extra info on $\mI$, and $\ell$}

Keeping the same notation as before, we recall that $[K,K] = [C,D] = \{ (cd)^2 \vert c \in C, d \in D\}$. Each $(cd)^2 = cdcd = [c,d]$ is an element in $\mS$
which obviously is an involution if it is nontrivial, so $[K,K] \leq \mI$. Also, note that the surjective morphism 
\begin{equation}
\omega:\ \mS \ \mapsto\ \mS^2:\ s \ \mapsto\ s^2,
\end{equation}
has kernel $\mI$, so that $\mS^2 \cong \mS/\mI$. It follows that $\vert [ K, K ] \vert = \vert \mS^2 \vert = t/\ell$. 
It now follows that
\begin{equation}
\label{kk}
\ell = \vert \mI \vert\ \geq\ \Big\vert [K,K] \Big\vert = \frac{t}{\ell},
\end{equation}
so that $\ell \geq \sqrt{t}$. \\

\section{Solution of Frohardt's problem}
\label{Frohthm}

Recall that we suppose that $\Gamma$ is not a TGQ. 
For this part of the paper we will invoke a result of Ott \cite{Ott}. On p. 413 of his paper (Proposition 4.2), Ott considers a maximal subgroup $M$ of $\Phi(K)$ (in his paper, the notation $G$ is used for our $K$ and $U_0$ is used for $\mS$). Then \cite[Corollary 4.5, p. 414]{Ott} says that $\vert M \vert = t/4$. In this section, we only will use that observation of his paper. We do want to obtain a full solution of Frohardt's problem, so we cannot rely on the main result of \cite{Ott}: there, the case that $t$ is a square is excluded, and we want to handle that too. Besides, after \cite[Corollary 4.5, p. 414]{Ott}, the proof of the main result still takes about 20 pages, and we can rely on equation (\ref{kk}) instead to obtain a very short approach to Frohardt's problem. 

So let $\vert M \vert = t/4$; then as a maximal subgroup of the $2$-group $\Phi(K)$, it has index 2, so that $\vert \Phi(K) \vert = t/2$. But as $\Phi(K) = [K, K]$, equation (\ref{kk}) tells us that $\vert \Phi(K) \vert = t/\ell$, so that $\ell = 2$. As $\sqrt{t} \leq \ell$ by the previous section, it follows that $t \leq 4$. For the cases $t = 2$ and $t = 4$, it is well known 
that only the classical quadrangles $\mW(2)$ and $\mW(4)$ occur (see \cite[chapter 6]{PT}), which means that $\Gamma$ is a TGQ. 

This shows that the class (3) in Frohardt's Theorem \ref{Fro88} is empty, as all hypothetical members are TGQs, implying that $K$ is elementary abelian. \eop \\

\section{STGQs of even order $t$, $t$ not a square, and cyclic STGQs}
\label{cyclic} 

{\bf In this section, we drop Frohardt's assumptions altogether, and initially work with general STGQs of order $(t, t)$.}  

In this section, we will also rely on one extra (nontrivial) result of Ott's paper \cite{Ott}, namely the fact that $\Phi(K) = \Phi(\mS)$ if $\Gamma$ is not a TGQ  (\cite[Theorem 4.11]{Ott}). It is our goal to give a very short alternative geometric approach to what follows after Theorem in Ott's paper (which is a long argument of 12 pages). Since we invoke that theorem, we need to assume that $t$ is not a square. Instead of ``directly'' obtaining the result, we choose a more interesting approach: we will completely classify cyclic STGQs of order $(t, t)$, cf. Theorem \ref{cyc} (without any restrictions on $t$).\\

One crucial observation is the following. 

\begin{lemma}
\label{cyclem}
We have that $\mS$ is cyclic. 
\end{lemma}

{\em Proof.}\quad
As the group $M$ of the previous section has size $t/4$ and is a maximal subgroup of $\Phi(K)$ (cf. Ott's paper), we have that $\vert \Phi(K) \vert = t/2$, and 
so $\vert \Phi(\mS) \vert = t/2$. 
As $[ \mS : \Phi(\mS)] = 2$, we have that $\mS$ is generated by one element (a finite $p$-group $P$ is cyclic if and only if the Frattini quotient $P/\Phi(P)$ is cyclic). \eop \\

\subsection*{Angles}

Let $(\Gamma, K)$ be an STGQ with elation point $x$, and associated group of symmetries (with center $x$) $\mS$. Let $(\mF,\mF^*)$ be the associated Kantor family, and let $A \in \mF$ be arbitrary. Let $A$ fix the line $U \sim [U] \I x$, where $U$ is not incident with $x$. Suppose that $A$ fixes each point incident with $[U]$. Define the {\em angle spectrum} of each element $a \in A^\times$ as follows. Let $V \sim [U]$, with $V$ not incident with $x$; then the {\em angle} of $a$ at $V$ is the unique element $s$ of $\mS$ for which $V^a = V^s$. Then the {\em angle spectrum} of $a$ is the set of all angles of $a$ at $V$, where $V$ varies over the lines 
in $[U]^\perp$ which are not incident with $x$.  \\



We now prove the following theorem. 

\begin{theorem}[Cyclic STGQs]
\label{cyc}
Suppose $(\Gamma,K)$ is an STGQ of order $t$, and suppose that $\mS$ is cyclic. Then $t = p$ is a prime and $\Gamma \cong \mW(p)$. 
\end{theorem}

{\em Proof.}\quad 
For $t$ odd, we already know that $\Gamma \cong \mW(t)$ and as $\mS$ is elementary abelian in this case, it follows that $t$ is prime (in which case $\mS \cong C_t$). 

So now suppose that $s = t$ is even, and put $\mS = \langle \alpha \rangle$. We suppose that $t \geq 4$. 
Let $y \not\sim x$ be some fixed point of $\Gamma$, and if $U \I y$, let $x_U := [U] \cap U$, where $[U]$ is the projection of $U$ on $x$. Define the elements of the Kantor family of $(\Gamma,K)$ as  $\mF = \{ K_U\ \vert\ U \I y \}$ and $\mF^* = \{ K_{x_U}\ \vert\ U \I y \}$. 
Let $A \in \mF$ be arbitrary, with $A = K_U$ and $U \I y$,  and let $A$ act on $\mS$ by conjugation. Define $N(A)$ to be the kernel of this action and note that each element of $N(A)$ commutes with $\mS$ and hence fixes $x_U$ linewise. Then $A/N(A)$ is a subgroup of the automorphism group of $\mS$, which is well known to be isomorphic to $C_2 \times C_{2^{h - 2}}$. But as $A/N(A)$ is elementary abelian, this implies that  
\begin{equation}
A/N(A) \leq C_2 \times C_2,
\end{equation} 
where the latter is the $2$-torsion subgroup of $C_2 \times C_{2^{h - 2}}$. It follows that $\vert N(A) \vert \in \{ t/4, t/2, t \}$. \\


Let $U, x_U, [U], A = K_U$ be as before, with $U$ not incident with $x$ arbitrary; then $N(A) \leq A$ has size $t, t/2$ or $t/4$. Now let $y \ne x_U$ be any other point on $[U]$ different from $x$. \\

\subsection*{Any element of $N(A)^\times$ fixes $t/2 + 1$ points on $[U]$ linewise, and all angles are given by $\id$ or $\alpha^{t/2}$}

Let $a \in N(A)^\times$. As $a$ fixes the point $x_U$ linewise, it commutes with every element in $\mS$. Suppose $\beta$ is an angle of $a$; then $\beta a$ fixes a line (not incident with $x$), so it is an involution. It follows that $(\beta a)^2 = \beta^2 = \id$. So $\beta$ is the unique involution in $\mS$, $\alpha^{t/2}$. 
Now suppose $a$ fixes $\lambda$ points on $[U] \setminus \{x\}$ linewise, so that $t - \lambda$ points remain in $[U] \setminus \{x\}$ which are all incident with precisely one fixed line (which is $[U]$). If we would have that $\lambda < t/2$, then $t/2 < t - \lambda < t$. Considering the element $a\alpha^{t/2}$, we see that 
it fixes $t - \lambda$ points in $[U] \setminus \{ x\}$ linewise. As we will observe in the next paragraph, this is a contradiction as $t - \lambda$ has to be a power of $2$. 
\hfill{\bf QED} \\

\subsection*{The map $\zeta$}

Now let $a \in N(A)^\times$ be arbitrary; then $a$ fixes $t/2$ points on $[U]$ (besides $x$) linewise. Call this point set $\mB(a)$. Let $B \ne A$ in $\mF$; then it is easy to see that $\vert C_B(a) \vert = t/2$, and the set $\mB(a)$ can be identified with a unique subgroup of $B$ (which acts sharply transitively on $\mB(a)$). As $B$ is an elementary abelian $2$-group (of order $t = 2^h$), we can associate an $(h - 1)$-dimensional projective space over $\F_2$ with $B$, which we denote by $\Pi(B)$; we then have that $\mB(a)$ defines 
a hyperplane $\Pi(a)$ in $\Pi(B)$. Consider the map
\begin{equation}
\zeta:\ N(A)\ \mapsto \Pi(B):\ a \ \mapsto\ \Pi(a), 
\end{equation} 
and note that $\zeta(\id) = \Pi(B)$. \\

\ul{For now, we suppose that this map is an injection} (if it is not, we will later easily show that a contradiction arises).

Note that the dual spaces of $\zeta(N(A)^\times)$ are contained in (in fact, form) a subspace of $\Pi(B)^* \simeq  \Pi(B)$ (where $\Pi(B)^*$ is the dual space of $\Pi(B)$). If $a$ and $b$ are in $N(A)^\times$, then $ab = c$ is such that $\zeta(c)$ is the third hyperplane on $\zeta(a) \cap \zeta(b)$; in the dual space $\Pi(B)^*$, $\zeta(c)^*$ is the third point of the line $\zeta(a)^*\zeta(b)^*$. It follows that ${\zeta(N(A)^\times)}^* \simeq \zeta(N(A)^\times)$ is a subspace of $\Pi(B)^* \simeq \Pi(B)$. \\

\subsection*{Subgroups of $N(A)$ which fix at least one arbitrary point $x' \ne x_U, x$ linewise have at least size $t/8$} 

Let $A$ be arbitrary in $\mF$ as before, and keep using the same notation. Let $x' \ne x_U$ be arbitrary in $[U] \setminus \{ x\}$. Then the subgroup of $A$ 
of elements which fix both $x_U$ and $x'$ linewise, has at least size $t/8$. This follows from properties of the map $\zeta$ as we will see below. Since the dual spaces of the elements in $\zeta(N(A)^\times)$ form a subspace
of $\Pi(B)^* \simeq \Pi(B)$, it follows that precisely one extra point $\widetilde{x}$ in $[U] \setminus \{x_U,x\}$ is fixed linewise by $N(A)$ if $\vert N(A) \vert = t/2$ (${\zeta(N(A)^\times)}^*$ then is a hyperplane), or precisely three extra points $\widetilde{x}_1, \widetilde{x}_2, \widetilde{x}_3$ when $\vert N(A) \vert = t/4$ (${\zeta(N(A)^\times)}^*$ then is a subspace of codimension $2$). Suppose that $x'$ is not one of the extra points.  
We distinguish three cases. Below, we put $t = 2^h$, so that $\Pi(B)$ is a projective $\F_2$-space of dimension $h - 1$. \\

\begin{itemize}
\item[(0)] \ul{$\vert N(A) \vert = t$}.\quad 
Since all angles are in $\{ \id, \alpha^{t/2}\}$, the number of elements in $N(A)^\times$ which fix $x_U$ and $x'$ linewise is at least $t/2$. (Let $V \I x'$, $V \ne [U]$; then $\vert N(A)_V \vert \leq 2$.)  
\item[(1)] \ul{$\vert N(A) \vert = t/2$}.\quad
In this case all elements of $\zeta(N(A)^\times)$ contain one fixed point $r$ in $\Pi(B)$ (and all hyperplanes in $\Pi(B)$ containing that point are contained in $\zeta(N(A)^\times)$). The number of elements in $N(A)^\times$ which fix $x_U$ and $x'$ linewise is the number 
of sets $\mB(a)$ ($a \in N(A)^\times$) which contain besides $r$ a second point $r'$, and this equals the number of hyperplanes in $\Pi(B)$ which contain two different points, hence a line: $2^{h - 2} - 1 = t/4 - 1$.    
\item[(2)] \ul{$\vert N(A) \vert = t/4$}.\quad 
In this case all elements of $\zeta(N(A)^\times)$ contain one fixed line $R$ in $\Pi(B)$ (and all hyperplanes in $\Pi(B)$ containing that line are contained in $\zeta(N(A)^\times)$). The number of elements in $N(A)^\times$ which fix $x_U$ and $x'$ linewise is the number 
of sets $\mB(a)$ ($a \in N(A)^\times$) which contain both $x_U$ and $x'$, and this equals the number of hyperplanes in $\Pi(B)$ which contain a given point and an extra line $R$, hence a plane: $2^{h - 3} - 1 = t/8 - 1$.   
\end{itemize}

\hfill{\bf QED} \\

\subsection*{Alternative proof using angles} 

We have seen that all elements in $N(A)$ only allow angles contained in $\{ \id, \alpha^{t/2} \}$. Take any line $Y$ on $x'$ different from $[U]$. Then $\vert N(A)_Y \vert \geq \vert N(A) \vert/2$, since at most one nontrivial angle can occur. And since $[N(A),\mS] = \{ \id \}$, it also follows that $N(A)_Y$ fixes all lines on $x'$. This gives us the same cases (0), (1), (2) as above.  \hfill{\bf QED} \\

In the rest of the proof, we suppose that $t \geq 16$ (if $t = 8$, we refer to \cite{Leug}).

\subsection*{Each element in $A$, where $A \in \mF$ is arbitrary, has an angle spectrum which is a subset of $\{ \id, \alpha^{t/2} \}$; as a consequence we have that $[k, a] \in \{ \id, \alpha^{t/2} \}$ for all $k \in K$, $a \in A$, $A \in \mF$}


First suppose that $a \in N(A)^\times$. Then we already know that only the angles $\id$ and $\alpha^{t/2}$ can occur. Now let $k \in K^\times$ be arbitrary, and consider $V^k$. If $a$ fixes $V^k$, we have that $[k,a] = \id$ since it is a symmetry with center $x$ which fixes $V^k$; 
if $a$ does not fix $V^k$, we have that $[k,a] = \alpha^{t/2}$, since it is a symmetry with center $x$ which sends $V^k$ to ${\Big(V^{k}\Big)}^{\alpha^{t/2}}$. So if $N(A) = A$, the statement is already proved.\\

Suppose that $a \in A^\times \setminus N(A)^\times$; then $a$ induces (by conjugation) a nontrivial involutory automorphism of $\mS$. Let $\gamma \ne \id$ be an angle of $a$. Then $\gamma a$ fixes a line not incident with $x$, so that it is an involution. 
Then $\gamma a \gamma a = \id$, so that $\gamma^a = \gamma^{-1}$. If $a$ and $\gamma$ commute, it follows that $\gamma$ is an involution, so $\gamma = \alpha^{t/2}$. Now suppose $a$ and $\gamma$ do not commute. Consider the line $U \I x_U$, $U \ne [U]$, which is supposed to be fixed by $a$. We identify 
the action of $a$ on $\mS$ with the action of $a$ on the lines incident with $x_U$ and different from $[U]$ (which we can do since $\mS$ acts 
sharply transitively on these lines: if $W \I x_U$ and $W \ne [U]$ with $W = V^\alpha$, then we identify $W$ with $\alpha \in \mS$). 
It is easy to calculate that the only nontrivial involutions in $\Aut(\mS)$ are given by the following images of $\alpha$:
\begin{quote}
\item
$\nu_1:\ \alpha:\ \mapsto\ \alpha^{t - 1}$; 
\item
$\nu_2:\ \alpha:\ \mapsto\ \alpha^{t/2 + 1}$; 
\item
$\nu_3:\ \alpha:\ \mapsto\ \alpha^{(t - 1)(t/2 + 1)} = \alpha^{t/2 - 1}$. 
\end{quote}
In the first case, $\nu_1$ fixes precisely two elements ($\id$ and $\alpha^{t/2}$) of $\mS$; in the second case $t/2$ elements (all $\alpha^n$ with $n$ even), and 
in the last case again two elements ($\id$ and $\alpha^{t/2}$). If $a$ is of the first type or the third, at most two lines different from $[U]$ are fixed per point 
in $[U] \setminus \{x\}$  (in fact, one can easily show that it is {\em precisely} two per such point, but we do not need this equality). If $a$ is of the second type, 
it fixes precisely $t/2$ lines incident with $x_U$ and different from $[U]$ (and many more on other points as well, as we will later see). 

Now suppose $a \in A^{\times}$ is of type 1 or type 3. Suppose $a$ fixes the line $U$ (not incident with $x$), and let $u^a = v$ with $u \I U$, $u$ not collinear with $x$. We suppose without loss of generality that $u$ ``defines the Kantor family'' | that is, if $W \I u$, then $K_W$ is an element of $\mF$ (and in particular, 
$K_U = A$). 
For any such $W \ne U$, if $B = K_W$ in $\mF$ and $[W]$ is the line incident with $x$ which meets $W$, define $U(B) = W$ as the line incident with $u$ which meets $[W]$, and let $V(B)$ be the line incident with $v$ which meets $[W]$. Let $u_B := U(B) \cap [W]$ and $v_B := V(B) \cap [W]$. 
Then as we have seen in a previous step, we have that $B_{[u_B,v_B]}$, the subgroup of $B$ which fixes $u_B$ and $v_B$ linewise, has size 
at least $t/8$. Note that $U^{B_{[u_B,v_B]}}$ consists of lines in $\{ [U], U(B), V(B) \}^{\perp}$. As a consequence, since $a$ maps 
$u$ to $v$ and fixes all points incident with $[U]$, it follows that $a$ fixes all lines in  $U^{B_{[u_B,v_B]}}$. Letting $B$ vary in $\mF \setminus \{A\}$, we 
see that $a$ has at least $t(t/8 - 1) + 1$ fixed lines which are not collinear with $x$. Supposing that $t \geq 32$ (the case $t = 16$ will be handled separately below), we see that type 1 and 3 are indeed not possible (since 
such elements would fix at most $2t$ such lines). 

Finally, suppose that $a$ is of type 2. Suppose $\gamma \ne \id$ is an angle of $a$; then we have seen that $\gamma^a = \gamma^{-1}$. If we inspect the specific form of $a$ and suppose that $\gamma = \alpha^n$, we get that 
\begin{equation}
\alpha^{n(t/2 + 1)} \ =\ \alpha^{-n},
\end{equation} 
so that $\alpha^{n(t/2 + 2)} = \id$. This is only possible if $n = t$ ($\gamma = \id$) or $n = t/2$ ($\gamma = \alpha^{t/2}$).
Now let $k \in K^\times$ be arbitrary, and consider $V^k$. Similarly as in the case $a \in N(A)^\times$, we have that either $[k,a] = \id$ or $[k,a] = \alpha^{t/2}$. 

The claim is proved. 
\hfill{\bf QED} \\

\subsection*{The case $t = 16$}

We use the notation of the preceding part of the proof. Suppose again that $a \in A^{\times}$ is of type 1 or type 3. 
Suppose that $\vert N(B) \vert = t/4$; then $B$ contains elements of every type. Suppose $\beta \in B$ is an element of type 2. We have seen moments ago that for such an element, only angles in $\{ \id, \alpha^{t/2} \}$ can arise. Whence the same is true for all elements in $\langle N(B), \beta \rangle$ (which is a group of size $2\vert N(B) \vert$). It follows that $\vert B_{U(B), V(B)} \vert \geq \vert \langle N(B), \beta \rangle \vert/2 = \vert N(B) \vert$. Note that $U^{B_{U(B),V(B)}}$ consists of lines in $\{ [U], U(B), V(B) \}^{\perp}$.
Hence in the counting argument above, we can actually replace $t/8$ by $t/4$ (such as in the cases $\vert N(B) \vert = t/2$ and $\vert N(B) \vert = t$), and we obtain that $a$ has at least $t(t/4 - 1) + 1$ fixed lines which are not collinear with $x$. When $t = 16$, this amount is strictly larger than $2t$, so we have a contradiction. When $a$ is of type 2, the argument is the same as above. So if $k \in K$, we have that either $[k,a] = \id$ or $[k,a] = \alpha^{t/2}$ for any $a \in A$. 
\hfill{\bf QED} \\

\subsection*{We have that $[K, K] \leq C_2$ (up to isomorphism)}

This is easy now: consider any commutator $[u, v]$ with $u, v \in K$; then $u$ and $v$ are products of elements contained in the members of $\mF$. It follows that such a commutator is a product of (conjugates of) elements $[e, f]$, with $e, f$ elements of members in $\mF$. Whence $[u, v]$ is a product of copies 
of $\alpha^{t/2}$. The claim follows. \hfill{\bf QED} \\

\subsection*{We have that $\Gamma \cong W(p)$ with $s = p$ a prime}

Let $\alpha \in \mS$. Then we can write $\alpha$ as $a_1a_2 \cdots a_k$ with each $a_k$ contained in some element of $\mF$; it follows that 
\begin{equation}
\alpha^2\ =\ a_1^2a_2^2\cdots a_k^2\ell\ =\ \ell,
\end{equation}
with $\ell$ a product of commutators, implying that the order of each such $\alpha$ is at most $4$. This means that $s = \vert \mS \vert \leq 4$ as $\mS$ is cyclic. 
By \cite{Leug} we conclude that $\Gamma \cong \mW(s)$ with $s \in \{ 2, 4\}$. Only in the case $\Gamma \cong \mW(2)$ we have that $\mS$ is cyclic. \hfill{\bf QED}\\

Now suppose that $\zeta$ is not an injection. Let $a, \widetilde{a}$ be different nontrivial elements of $N(A)$ for which $\zeta(a) = \zeta(\widetilde{a})$. Then as 
they have the same spectrum, it follows that $a\widetilde{a}$ is a symmetry with axis $[U]$. Let $S$ be the ``kernel'' of $\zeta$: this is the subgroup of symmetries of $N(A)$ with axis $[U]$. Now $\Big(N(A)/S \Big)^\times$ induces a subspace of $\Pi(B)^*$ (instead of $N(A)^\times$). Putting $\vert N(A) \vert = 2^{h - m}$, $m \in \{ 0, 1, 2  \}$ and $\vert S \vert = 2^k$, and using the same notation as before, we obtain that $2^{h - m - k - 1} - 1$ elements of $\Big(N(A)/S \Big)^\times$ fix $x_U$ and $x'$ linewise (note that elements of $N(A)/S$ act naturally on the lines in $[U]^\perp$). As each element of $S$ also fixes $x_U$ and $x'$ linewise, we obtain that $(2^{h - m - k - 1} - 1)2^k + 2^k - 1 = 2^{h - m - 1} - 1$ nontrivial elements of $N(A)$ fix $x_U$ and $x'$ linewise, which is the exact same result as in the injective case. Note that the short alternative approach using angles also works here!
Now the rest of the proof is the same. It follows that $\Gamma \cong \mW(2)$. \\


The proof is complete. 

\eop \\

Theorem \ref{cyc} immediately leads to a very short proof of Ott's result:

\begin{corollary}[Ott's Result]
Let $(\Gamma,K)$ be an STGQ of even order $t$. If $t$ is not a square, then $\Gamma$ is a TGQ and $K$ is elementary abelian. 
\end{corollary}

{\em Proof}.\quad
Suppose $\Gamma$ is not a TGQ, so that $K$ is not elementary abelian. 
Then by Lemma \ref{cyclem} we know that $\mS$ is cyclic. So by Theorem \ref{cyc} we conclude that $\Gamma \cong \mW(2)$, contradiction (since that {\em is} a TGQ). \eop \\

\end{document}